\documentclass[12pt,a4paper]{article}
\usepackage{amsmath}
\usepackage{amssymb}
\usepackage{t1enc}
\usepackage[latin1]{inputenc}
\usepackage[english]{babel}
\usepackage{amsfonts}
\usepackage{latexsym}
\usepackage{bm}
\newtheorem{theorem}{Theorem}[section]
\newtheorem{prop}[theorem]{Proposition}
\newtheorem{lemma}[theorem]{Lemma}
\newtheorem{cor}[theorem]{Corollary}

\newtheorem{rem}[theorem]{Remark}
\newtheorem{eg}[theorem]{Example}
\newtheorem{construction}[theorem]{Construction}

\begin{document}
\title{The maximum sum and maximum product of sizes of
cross-intersecting families}

\author{Peter Borg\\[5mm]
Department of Mathematics, University of Malta, Malta\\
\texttt{p.borg.02@cantab.net}}
\date{December 29, 2012}
\maketitle
%\date{13th June 2010}

\begin{abstract}
We say that a set $A$ \emph{$t$-intersects} a set $B$ if $A$ and
$B$ have at least $t$ common elements. A family $\mathcal{A}$ of
sets is said to be \emph{$t$-intersecting} if each set in
$\mathcal{A}$ $t$-intersects any other set in $\mathcal{A}$.
Families $\mathcal{A}_1, \mathcal{A}_2, ..., \mathcal{A}_k$ are
said to be \emph{cross-$t$-intersecting} if for any $i$ and $j$
in $\{1, 2, ..., k\}$ with $i \neq j$, any set in $\mathcal{A}_i$
$t$-intersects any set in $\mathcal{A}_j$. We prove that for any
finite family $\mathcal{F}$ that has at least one set of size at
least $t$, there exists an integer $\kappa \leq |\mathcal{F}|$
such that for any $k \geq \kappa$, both the sum and the product
of sizes of any $k$ cross-$t$-intersecting sub-families
$\mathcal{A}_1, ..., \mathcal{A}_k$ (not necessarily distinct or
non-empty) of $\mathcal{F}$ are maxima if $\mathcal{A}_1 = ... =
\mathcal{A}_k = \mathcal{L}$ for some largest $t$-intersecting
sub-family $\mathcal{L}$ of $\mathcal{F}$. We then study the
smallest possible value of $\kappa$ and investigate the case $k <
\kappa$; this includes a cross-intersection result for straight lines that demonstrates that it is possible to have $\mathcal{F}$ and $\kappa$ such that for any $k < \kappa$, the configuration $\mathcal{A}_1 = ...
= \mathcal{A}_k = \mathcal{L}$ is neither optimal for the sum nor
optimal for the product. We also outline solutions for various
important families $\mathcal{F}$, and we provide solutions for
the case when $\mathcal{F}$ is a power set.
\end{abstract}

\section{Introduction}\label{Intro}
%Give the parameter $\alpha_{\mathcal{F}}$ some suitable name.
%Maybe "cross-intersection parameter"?

Unless otherwise stated, throughout this paper we shall use small
letters such as $x$ to denote elements of a set or positive
integers, capital letters such as $X$ to denote sets, and
calligraphic letters such as $\mathcal{F}$ to denote
\emph{families} (that is, sets whose elements are sets
themselves). Unless specified, sets and families are taken to
be finite and may be the \emph{empty set} $\emptyset$. An \emph{$r$-set} is a set of size $r$, that is, a set having
exactly $r$ elements. For any integer $n \geq 1$, $[n]$ denotes
the set $\{1, ..., n\}$ of the first $n$ positive integers.

Given an integer $t \geq 1$, we say that a set $A$
\emph{$t$-intersects} a set $B$ if $A$ and $B$ have at least $t$
common elements. A family $\mathcal{A}$ is said to be
\emph{$t$-intersecting} if each set in $\mathcal{A}$
$t$-intersects any other set in $\mathcal{A}$ (i.e.~$|A \cap B|
\geq t$ for any $A, B \in \mathcal{A}$ with $A \neq B$). A
$1$-intersecting family is also simply called an
\emph{intersecting family}. Families $\mathcal{A}_1, ...,
\mathcal{A}_k$ are said to be \emph{cross-$t$-intersecting} if
for any $i$ and $j$ in $[k]$ with $i \neq j$, any set in $\mathcal{A}_i$ $t$-intersects any set in $\mathcal{A}_j$ (i.e.~$|A \cap B| \geq t$ for any $A \in \mathcal{A}_i$ and any $B \in \mathcal{A}_j$). Cross-$1$-intersecting families are also simply called
\emph{cross-intersecting families}.

Let ${[n] \choose r}$ denote the family of all subsets of
$[n]$ of size $r$. The classical Erd\H os-Ko-Rado (EKR) Theorem \cite{EKR} says that if $n$ is sufficiently larger than $r$, then a
$t$-intersecting sub-family $\mathcal{A}$ of ${[n] \choose r}$
has size at most ${n-t \choose r-t}$, which is the number of sets
in the $t$-intersecting sub-family of ${[n] \choose r}$
consisting of those sets having $[t]$ as a subset. The EKR Theorem
inspired a wealth of results of this kind, that is, results that
establish how large a system of sets can be under certain
intersection conditions; the survey papers \cite{DF,F} are
recommended.

For $t$-intersecting sub-families of a given family $\mathcal{F}$,
the natural question to ask is how large they can be. For
cross-$t$-intersecting families, two natural parameters arise: the
sum and the product of sizes of the cross-$t$-intersecting
families (note that the product of sizes of $k$ families
$\mathcal{A}_1, ..., \mathcal{A}_k$ is the number of $k$-tuples
$(A_1, ..., A_k)$ such that $A_i \in \mathcal{A}_i$ for each $i
\in [k]$). It is therefore natural to consider the problem of
maximising the sum or the product of sizes of $k$
cross-$t$-intersecting sub-families (not necessarily distinct or
non-empty) of a given family $\mathcal{F}$.

The main result in this paper (Theorem~\ref{main} below) relates
both the maximum sum and the maximum product of sizes of $k$
cross-$t$-intersecting sub-families of a family $\mathcal{F}$ to
the maximum size of a $t$-intersecting sub-family of $\mathcal{F}$
when $k$ is not smaller than a certain value depending on
$\mathcal{F}$ and $t$. It gives the maximum sum and the maximum
product in terms of the size of a largest $t$-intersecting
sub-family.

For any non-empty family $\mathcal{F}$, let $\alpha(\mathcal{F})$
denote the size of a largest set in
$\mathcal{F}$. Suppose $\alpha(\mathcal{F}) < t$, and let $\mathcal{A}_1, ...,
\mathcal{A}_k$ ($k \geq 2$) be sub-families of $\mathcal{F}$.
Then $\mathcal{A}_1, ..., \mathcal{A}_k$ are
cross-$t$-intersecting if and only if at most one of them is
non-empty (since no set in $\mathcal{F}$ $t$-intersects itself or
another set in $\mathcal{F}$). Thus, if
$\mathcal{A}_1, ..., \mathcal{A}_k$ are cross-$t$-intersecting,
then the product of their sizes is $0$ and the sum of their sizes
is at most the size $|\mathcal{F}|$ of $\mathcal{F}$ (which is
attained if and only if one of them is $\mathcal{F}$ and the
others are all empty). This completely solves our problem for the
case $\alpha(\mathcal{F}) < t$.

We now address the case $\alpha(\mathcal{F}) \geq t$. Before stating our main result, we need to introduce some definitions and parameters.

For any family $\mathcal{A}$, let $\mathcal{A}^{t,+}$ be the
($t$-intersecting) sub-family of $\mathcal{A}$ given by
\[\mathcal{A}^{t,+} = \{A \in
\mathcal{A} \colon |A \cap B| \geq t \mbox{ for any } B \in
\mathcal{A} \mbox{ with } A \neq B\},\]
and let
\[\mathcal{A}^{t,-} = \mathcal{A} \backslash \mathcal{A}^{t,+}.\]
In simple terms, a set $A$ in $\mathcal{A}$ is in
$\mathcal{A}^{t,-}$ if there exists a set $B$
in $\mathcal{A}$ such that $A \neq B$ and $A$ does not $t$-intersect $B$, otherwise $A$ is in $\mathcal{A}^{t,+}$. The
definitions of $\mathcal{A}^{t,+}$ and $\mathcal{A}^{t,-}$ are
generalisations of the definitions of $\mathcal{A}^*$ and
$\mathcal{A}'$ in \cite{Borg4, Borg3, Borg2, Borg5, BL2}; $\mathcal{A}^*
= \mathcal{A}^{1,+}$ and $\mathcal{A}' = \mathcal{A}^{1,-}$.

Let $l(\mathcal{F},t)$ denote the size of a largest
$t$-intersecting sub-family of a non-empty family
$\mathcal{F}$. For any sub-family $\mathcal{A}$ of $\mathcal{F}$,
we define
\[ \beta(\mathcal{F},t,\mathcal{A}) = \begin{cases}
\displaystyle \frac{l(\mathcal{F},t) -
|\mathcal{A}^{t,+}|}{|\mathcal{A}^{t,-}|} & \text{if
$\mathcal{A}^{t,-} \neq \emptyset$;}\\[4mm]
\displaystyle \frac{l(\mathcal{F},t)}{|\mathcal{F}|} & \text{if
$\mathcal{A}^{t,-} = \emptyset$;}
\end{cases} \]
so $|\mathcal{A}^{t,+}| + \beta(\mathcal{F},t,\mathcal{A})
|\mathcal{A}^{t,-}| \leq l(\mathcal{F},t)$ (even if
$\mathcal{A}^{t,-} = \emptyset$, because $|\mathcal{A}^{t,+}| \leq
l(\mathcal{F},t)$ since $\mathcal{A}^{t,+}$ is $t$-intersecting).
We now define
\[\beta(\mathcal{F}, t) = \min \{\beta(\mathcal{F}, t,
\mathcal{A}) \colon \mathcal{A} \subseteq \mathcal{F}\}.\]
Therefore,
\begin{equation} |\mathcal{A}^{t,+}| + \beta(\mathcal{F},t)
|\mathcal{A}^{t,-}| \leq l(\mathcal{F},t) \quad \mbox{for any
$\mathcal{A} \subseteq \mathcal{F}$.} \label{label1}
\end{equation}
In Section~\ref{betasection} we show that in fact
\begin{equation} \beta(\mathcal{F},t) = \max\left\{c \in
\mathbb{R} \colon c \leq \frac{l(\mathcal{F},t)}{|\mathcal{F}|},
\, |\mathcal{A}^{t,+}| + c|\mathcal{A}^{t,-}| \leq
l(\mathcal{F},t) \mbox{ for any } \mathcal{A} \subseteq
\mathcal{F}\right\} \label{label1.5}
\end{equation}
(see Proposition~\ref{beta1}), where $\mathbb{R}$ is the set of
real numbers, and we also determine other basic facts about the
parameter $\beta(\mathcal{F},t)$; in particular, we show that we
actually have
\begin{equation} \frac{1}{|\mathcal{F}|} \leq \beta(\mathcal{F},t)
\leq \frac{l(\mathcal{F},t)}{|\mathcal{F}|} \label{label2}
\end{equation}
(see Propositions~\ref{beta2} and \ref{beta1}). In
Section~\ref{betasection2} we point out various important
families $\mathcal{F}$ for which $\beta(\mathcal{F},t)$ is known
to be $\frac{l(\mathcal{F},t)}{|\mathcal{F}|}$.

By the lower bound in (\ref{label2}), for any non-empty family
$\mathcal{F}$, we can define
\[\kappa(\mathcal{F},t) = \frac{1}{\beta(\mathcal{F},t)}\] and we
have
\[\kappa(\mathcal{F},t) \leq |\mathcal{F}|.\]
We can now state our main result.

\begin{theorem} \label{main} Let $\mathcal{A}_1, ..., \mathcal{A}_k$ be cross-$t$-intersecting sub-families of a family $\mathcal{F}$ with $\alpha(\mathcal{F}) \geq t$. If $k \geq \kappa(\mathcal{F},t)$, then
\[\sum_{i = 1}^k |\mathcal{A}_i| \leq k(l(\mathcal{F},t)) \quad
\mbox{ and } \quad \prod_{i = 1}^k |\mathcal{A}_i| \leq \left(
{l(\mathcal{F},t)} \right)^k,\]
and both bounds are attained if $\mathcal{A}_1 = ... =
\mathcal{A}_k = \mathcal{L}$ for some largest $t$-intersecting
sub-family $\mathcal{L}$ of $\mathcal{F}$. Moreover, if $k >
\kappa(\mathcal{F},t)$, then in both inequalities, equality holds only if $\mathcal{A}_1 = ... = \mathcal{A}_k = \mathcal{L}$ for some largest $t$-intersecting sub-family $\mathcal{L}$ of
$\mathcal{F}$.
\end{theorem}

In Section~\ref{sumproblem} we prove the following result, which
tells us that if $k < \kappa(\mathcal{F},t)$, then the sum
inequality above does not hold for $\mathcal{A}_1, ...,
\mathcal{A}_k$ with a maximum value of $\sum_{i = 1}^k
|\mathcal{A}_i|$.
%(see Theorem~\ref{summax2}).

\begin{theorem}\label{summax2} Let $\mathcal{F}$ be a family with
$\alpha(\mathcal{F}) \geq t$. Let $\mathcal{A}_1, ...,
\mathcal{A}_k$ be cross-$t$-intersecting sub-families of
$\mathcal{F}$ such that $\sum_{i=1}^k |\mathcal{A}_i|$ is maximum. Then:\\
(i) $\sum_{i=1}^k |\mathcal{A}_i| = k(l(\mathcal{F},t))$ if $k \geq
\kappa(\mathcal{F},t)$; \\
(ii) $\sum_{i=1}^k |\mathcal{A}_i| > k(l(\mathcal{F},t))$ if $k <
\kappa(\mathcal{F},t)$.
\end{theorem}

In Theorem~\ref{main} the product inequality follows immediately from the sum inequality and the following elementary result,
known as the Arithmetic Mean-Geometric Mean (AM-GM) Inequality.

\begin{lemma}[AM-GM Inequality] \label{AMGM} If $x_1, x_2, \dots,
x_k$ are non-negative real numbers, then
\[\left( \prod_{i=1}^k x_i \right)^{1/k} \leq \frac{1}{k}\sum_{i=1}^k x_i.\]
\end{lemma}
Indeed, if $\sum_{i=1}^k |\mathcal{A}_i| \leq k(l(\mathcal{F},t))$, then, by Lemma~\ref{AMGM}, $\left( \prod_{i=1}^k|\mathcal{A}_i|
\right)^{1/k} \leq l(\mathcal{F},t)$ and hence $\prod_{i=1}^k|\mathcal{A}_i| \leq (l(\mathcal{F},t))^k$. Therefore, if the configuration $\mathcal{A}_1 = ... =
\mathcal{A}_k = \mathcal{L}$ (where $\mathcal{L}$ is as in
Theorem~\ref{main}) gives a maximum sum, then $\mathcal{A}_1 =
... = \mathcal{A}_k = \mathcal{L}$ also gives a maximum product.
The converse is not true; indeed, as demonstrated in Section~\ref{productproblem}, we may have that $2 \leq k <
\kappa(\mathcal{F},t)$ and $\mathcal{A}_1 = ... = \mathcal{A}_k =
\mathcal{L}$ still gives a maximum product, in which case $\mathcal{A}_1 = ... = \mathcal{A}_h =
\mathcal{L}$ gives a maximum product for any $h \geq k$ (see Lemma~\ref{prodext}). However, in Section~\ref{productproblem} we prove the following interesting fact.

\begin{rem} \label{productthreshold} \emph{Just like the threshold $\kappa(\mathcal{F},t)$ for the maximum sum part of Theorem~\ref{main} cannot be improved (by Theorem~\ref{summax2}), the threshold $\kappa(\mathcal{F},t)$ for the maximum product part of Theorem~\ref{main} can neither be
improved in general. Indeed, we will give a (geometrical) construction of a family $\mathcal{F}$ such that for any $k <
\kappa(\mathcal{F},t)$, the product of sizes of $k$ cross-$t$-intersecting
sub-families $\mathcal{A}_1, ..., \mathcal{A}_k$ of $\mathcal{F}$
is not maximum if $\mathcal{A}_1 = ... = \mathcal{A}_k =
\mathcal{L}$; see Construction~\ref{geomcon} and
Theorem~\ref{geomthm}.}
\end{rem}

We conclude this section by mentioning that in
Sections~\ref{sumproblem} and \ref{productproblem} we provide
various general results about the maximum sum and the maximum
product, respectively, and we also outline solutions for various
important families.

\section{Proof of the main result}

The proof of Theorem~\ref{main} relies on the following lemma.

\begin{lemma} \label{intxint} Let $\mathcal{A}_1, ...,
\mathcal{A}_k$ be cross-$t$-intersecting families, and let
$\mathcal{A} = \bigcup_{i=1}^k \mathcal{A}_i$. Then \\
(i) $\mathcal{A}^{t,+} = \bigcup_{i=1}^k {\mathcal{A}_i}^{t,+}$,\\
(ii) $\mathcal{A}^{t,-} = \bigcup_{i=1}^k
{\mathcal{A}_i}^{t,-}$,\\
(iii) $|\mathcal{A}^{t,-}| = \sum_{i=1}^k
|{\mathcal{A}_i}^{t,-}|$.
\end{lemma}
\textbf{Proof.} Clearly $\mathcal{A}^{t,+} \subseteq
\bigcup_{i=1}^k {\mathcal{A}_i}^{t,+}$. Suppose $A \in
\bigcup_{i=1}^k {\mathcal{A}_i}^{t,+}$. Then $A \in
{\mathcal{A}_h}^{t,+}$ for some $h \in [k]$, meaning that $|A
\cap B| \geq t$ for any $B \in \mathcal{A}_h \backslash \{A\}$.
Also, by the cross-$t$-intersection condition, for any $j \in [k]
\backslash \{h\}$, $|A \cap B| \geq t$ for any $B \in
\mathcal{A}_j$. So $A \in \mathcal{A}^{t,+}$. Therefore,
$\bigcup_{i=1}^k {\mathcal{A}_i}^{t,+} \subseteq
\mathcal{A}^{t,+}$. Together with $\mathcal{A}^{t,+} \subseteq
\bigcup_{i=1}^k {\mathcal{A}_i}^{t,+}$, this gives us (i).

Clearly $\bigcup_{i=1}^k {\mathcal{A}_i}^{t,-} \subseteq
\mathcal{A}^{t,-}$. Suppose $A \in {\mathcal{A}}^{t,-}$. Then $A
\in \mathcal{A}_h$ for some $h \in [k]$, and $|A \cap B| < t$ for
some $B \in \mathcal{A} \backslash \{A\}$. By the
cross-$t$-intersection condition, $B \notin \mathcal{A}_j$ for
each $j \in [k] \backslash \{h\}$. So $B \in \mathcal{A}_h$ and
hence $A \in {\mathcal{A}_h}^{t,-}$. Therefore,
${\mathcal{A}}^{t,-} \subseteq \bigcup_{i=1}^k
{\mathcal{A}_i}^{t,-}$. Together with $\bigcup_{i=1}^k
{\mathcal{A}_i}^{t,-} \subseteq \mathcal{A}^{t,-}$, this gives us
(ii).

Suppose ${\mathcal{A}_i}^{t,-} \cap {\mathcal{A}_j}^{t,-} \neq
\emptyset$ for some $i$ and $j$ in $[k]$ with $i \neq j$. Let $A \in
{\mathcal{A}_i}^{t,-} \cap {\mathcal{A}_j}^{t,-}$. Having $A \in
{\mathcal{A}_i}^{t,-}$ means that there exists a set $B$ in
${\mathcal{A}_i}^{t,-}$ such that $|A \cap B| < t$; however,
since $A \in \mathcal{A}_j$, this contradicts the
cross-$t$-intersection condition. So ${\mathcal{A}_i}^{t,-} \cap
{\mathcal{A}_j}^{t,-} = \emptyset$ for any $i$ and $j$
in $[k]$ with $i \neq j$, meaning that ${\mathcal{A}_1}^{t,-}, ...,
{\mathcal{A}_k}^{t,-}$ are disjoint. Together with (ii), this
gives us (iii).~\hfill{$\Box$}\\
\\
\textbf{Proof of Theorem~\ref{main}.} Suppose $k \geq \kappa(\mathcal{F},t)$. Then $\beta(\mathcal{F},t)
\geq 1/k$. Let $\mathcal{A}$ be the sub-family of $\mathcal{F}$
given by the union $\bigcup_{i = 1}^k \mathcal{A}_i$. We have
\begin{align} \sum_{i=1}^k |\mathcal{A}_i| &=
\sum_{i=1}^k |{\mathcal{A}_i}^{t,-}| + \sum_{i=1}^k
|{\mathcal{A}_i}^{t,+}| \nonumber \\
&\leq |\mathcal{A}^{t,-}| + k|\mathcal{A}^{t,+}| \quad \quad
\mbox{(by Lemma~\ref{intxint})} \nonumber \\
&= k\left(|\mathcal{A}^{t,+}| +
\frac{1}{k}|\mathcal{A}^{t,-}|\right) \nonumber \\
&\leq k\left(|\mathcal{A}^{t,+}| + \beta(\mathcal{F},t)
|\mathcal{A}^{t,-}|\right) \nonumber \\
&\leq k(l(\mathcal{F},t)) \quad \quad (\mbox{by (\ref{label1})}) \label{main1}
\end{align}
and, by Lemma~\ref{AMGM} and (\ref{main1}),
\begin{equation} \prod_{i=1}^k |\mathcal{A}_i| \leq \left( \frac{1}{k}\sum_{i=1}^k |\mathcal{A}_i| \right)^k \leq \left( {l(\mathcal{F},t)} \right)^k. \label{main2}
\end{equation}

If ${\mathcal{A}_1} = ... = {\mathcal{A}_k} = \mathcal{L}$ for some largest $t$-intersecting sub-family $\mathcal{L}$ of $\mathcal{F}$, then obviously ${\mathcal{A}_1}, ..., {\mathcal{A}_k}$ are cross-$t$-intersecting, $\sum_{i=1}^k
|\mathcal{A}_i| = k(l(\mathcal{F},t))$ and $\prod_{i=1}^k
|\mathcal{A}_i| = (l(\mathcal{F},t))^k$. Now suppose $k > \kappa(\mathcal{F},t)$ and either $\sum_{i=1}^k |\mathcal{A}_i| = k(l(\mathcal{F},t))$ or $\prod_{i=1}^k |\mathcal{A}_i| = (l(\mathcal{F},t))^k$. If $\prod_{i=1}^k |\mathcal{A}_i| = (l(\mathcal{F},t))^k$, then $\sum_{i=1}^k |\mathcal{A}_i| = k(l(\mathcal{F},t))$ by (\ref{main2}). So $\sum_{i=1}^k |\mathcal{A}_i| = k(l(\mathcal{F},t))$. Thus in (\ref{main1}) we have equality throughout. It follows that $|\mathcal{A}^{t,-}| = 0$ (since $k > \kappa(\mathcal{F},t)$ implies that $\frac{1}{k} < \beta(\mathcal{F},t)$), and hence $\mathcal{A} = \mathcal{A}^{t,+}$. So $\mathcal{A}$ is a $t$-intersecting sub-family of $\mathcal{F}$. Since $\sum_{i=1}^k |\mathcal{A}_i| = k(l(\mathcal{F},t))$ and $\mathcal{A}_i \subseteq \mathcal{A}$ for each $i \in [k]$, it follows that $\mathcal{A}_1 = ... = \mathcal{A}_k = \mathcal{A}$ and $\mathcal{A}$ is a largest $t$-intersecting sub-family of $\mathcal{F}$.~\hfill{$\Box$}

\section{The parameter $\beta(\mathcal{F},t)$}
\label{betasection}

Theorems~\ref{main} and \ref{summax2} tell us that
$\lceil \kappa(\mathcal{F},t) \rceil$ is the smallest integer $k_0$ such that
for any $k \geq k_0$, the configuration $\mathcal{A}_1 = ... =
\mathcal{A}_k = \mathcal{L}$ (as in the theorems) is optimal for
the maximisation of both the sum and the product of sizes. So
$\kappa(\mathcal{F},t)$ is an important parameter and hence worth
investigating. But $\kappa(\mathcal{F},t)$ is simply defined to
be the reciprocal of $\beta(\mathcal{F},t)$, and hence we may
instead focus on $\beta(\mathcal{F},t)$.

In this section we first establish some basic facts on
$\beta(\mathcal{F},t)$ and then we provide the value of
$\beta(\mathcal{F},t)$ for various important families
$\mathcal{F}$.

\subsection{General facts}

We start by proving the lower bound in (\ref{label2}) and characterising the families for which it is attained.

\begin{prop} \label{beta2} For any family $\mathcal{F} \neq
\emptyset$,
%and any %positive integer $t$,
%
\begin{equation} \beta(\mathcal{F},t) \geq
\frac{1}{|\mathcal{F}|}, \nonumber
\end{equation}
and equality holds if and only if $|A \cap B| < t$ for any
distinct $A$ and $B$ in $\mathcal{F}$.
\end{prop}
\textbf{Proof.} Since $\mathcal{F}$ is non-empty,
$l(\mathcal{F},t) \geq 1$ because any sub-family of $\mathcal{F}$
consisting of only one set is $t$-intersecting (by definition).
Let $\mathcal{A} \subseteq \mathcal{F}$. If $\mathcal{A}^{t,-} =
\emptyset$, then $\beta(\mathcal{F},t,\mathcal{A}) =
\frac{l(\mathcal{F},t)}{|\mathcal{F}|} \geq
\frac{1}{|\mathcal{F}|}$, and equality holds only if
$l(\mathcal{F},t) = 1$. Now suppose $\mathcal{A}^{t,-} \neq
\emptyset$. Let $A$ be a set in $\mathcal{A}^{t,-}$. Then
$\mathcal{A}^{t,+} \cup \{A\}$ is a $t$-intersecting sub-family
of $\mathcal{F}$, and hence $|\mathcal{A}^{t,+} \cup \{A\}| \leq
l(\mathcal{F},t)$. So $l(\mathcal{F},t) \geq |\mathcal{A}^{t,+}|
+ 1$. We therefore have $\beta(\mathcal{F},t,\mathcal{A}) =
\frac{l(\mathcal{F},t) -
|\mathcal{A}^{t,+}|}{|\mathcal{A}^{t,-}|} \geq
\frac{1}{|\mathcal{F}|}$, and equality holds only if
$l(\mathcal{F},t) - |\mathcal{A}^{t,+}| = 1$ and
$\mathcal{A}^{t,-} = \mathcal{F}$, in which case
$\mathcal{A}^{t,+} = \emptyset$ and hence $l(\mathcal{F},t) = 1$.

Therefore, $\beta(\mathcal{F},t) \geq \frac{1}{|\mathcal{F}|}$,
and equality holds only if $l(\mathcal{F},t) = 1$. Now clearly
$l(\mathcal{F},t) = 1$ if and only if $|A \cap B| < t$ for any
distinct $A$ and $B$ in $\mathcal{F}$, in which case either
$|\mathcal{F}| = 1$ or for all $\mathcal{A} \subseteq
\mathcal{F}$ with $|\mathcal{A}| \geq 2$, $\mathcal{A}^{t,-} =
\mathcal{A}$ and $\mathcal{A}^{t,+} = \emptyset$. So
$l(\mathcal{F},t) = 1$ implies that $\beta(\mathcal{F},t) =
\beta(\mathcal{F},t,\mathcal{F}) =
\frac{1}{|\mathcal{F}|}$.~\hfill{$\Box$}\\

We next prove (\ref{label1.5}), which gives a clear description of
$\beta(\mathcal{F}, t)$, and hence establish the upper bound in
(\ref{label2}).

\begin{prop} \label{beta1} For any family $\mathcal{F} \neq
\emptyset$, $\beta(\mathcal{F}, t)$ is the largest real number $c
\leq \frac{l(\mathcal{F},t)}{|\mathcal{F}|}$ such that
\[|\mathcal{A}^{t,+}| + c|\mathcal{A}^{t,-}|
\leq l(\mathcal{F},t) \quad \mbox{for any $\mathcal{A} \subseteq
\mathcal{F}$.}\]
\end{prop}
\textbf{Proof.} For any $\mathcal{A} \subseteq \mathcal{F}$,
$|\mathcal{A}^{t,+}| + \beta(\mathcal{F},t)|\mathcal{A}^{t,-}|
\leq |\mathcal{A}^{t,+}| +
\beta(\mathcal{F},t,\mathcal{A})|\mathcal{A}^{t,-}| \leq
l(\mathcal{F},t)$. Since $\beta(\mathcal{F},t,\emptyset) =
\frac{l(\mathcal{F},t)}{|\mathcal{F}|}$,
\begin{equation} \beta(\mathcal{F},t) \leq
\frac{l(\mathcal{F},t)}{|\mathcal{F}|}. \nonumber
\end{equation}
Suppose $\beta(\mathcal{F},t) <
\frac{l(\mathcal{F},t)}{|\mathcal{F}|}$. Let $d$ be a real
number greater than $\beta(\mathcal{F},t)$. Let ${\mathcal{A}_0}
\subseteq \mathcal{F}$ such that
$\beta(\mathcal{F},t,{\mathcal{A}_0}) = \beta(\mathcal{F},t)$.
Since $\beta(\mathcal{F},t,\mathcal{A}_0) \neq
\frac{l(\mathcal{F},t)}{|\mathcal{F}|}$, we have
${\mathcal{A}_0}^{t,-} \neq \emptyset$ and hence
$|{\mathcal{A}_0}^{t,+}| +
\beta(\mathcal{F},t,\mathcal{A}_0)|{\mathcal{A}_0}^{t,-}| =
l(\mathcal{F},t)$. So $|{\mathcal{A}_0}^{t,+}| +
d|{\mathcal{A}_0}^{t,-}| > l(\mathcal{F},t)$. Hence the
result.~\hfill{$\Box$}

\begin{rem} \emph{When $\mathcal{A}^{t,-} = \emptyset$ it does not
matter what $\beta(\mathcal{F},t,\mathcal{A})$ is, because
$|\mathcal{A}^{t,+}| + \beta(\mathcal{F},t,\mathcal{A})
|\mathcal{A}^{t,-}| = |\mathcal{A}^{t,+}| \leq l(\mathcal{F},t)$.
Thus we could define $\beta(\mathcal{F},t)$ to be the minimum
value of $\beta(\mathcal{F},t,\mathcal{A})$ such that
$\mathcal{A}^{t,-} \neq \emptyset$ when such a sub-family
$\mathcal{A}$ exists, i.e.~when $\mathcal{F}^{t,-} \neq \emptyset$
(i.e.~when $\mathcal{F}$ is not $t$-intersecting). However, this would still give us $\beta(\mathcal{F},t) \leq
\frac{l(\mathcal{F},t)}{|\mathcal{F}|}$; indeed,
\begin{align} \mathcal{F}^{t,-} \neq \emptyset \; &\Rightarrow \;
\beta(\mathcal{F},t,\mathcal{F}) = \frac{l(\mathcal{F},t) -
|\mathcal{F}^{t,+}|}{|\mathcal{F}^{t,-}|} =
\frac{l(\mathcal{F},t) - |\mathcal{F}^{t,+}|}{|\mathcal{F}| -
|\mathcal{F}^{t,+}|} \leq \frac{l(\mathcal{F},t)}{|\mathcal{F}|}
\label{label3} \\
&\Rightarrow \; \beta(\mathcal{F},t) \leq
\frac{l(\mathcal{F},t)}{|\mathcal{F}|}. \nonumber
\end{align}
}
\end{rem}

If $\mathcal{F}^{t,-} \neq \emptyset$ and $\mathcal{F}^{t,+} \neq
\emptyset$, then the inequality in (\ref{label3}) is strict. Thus,
if $\beta(\mathcal{F},t) = \frac{l(\mathcal{F},t)}{|\mathcal{F}|}$, then either $\mathcal{F}^{t,-} = \emptyset$ or $\mathcal{F}^{t,+} =\emptyset$. Therefore,
\begin{equation} \beta(\mathcal{F},t) =
\frac{l(\mathcal{F},t)}{|\mathcal{F}|} \; \Rightarrow \;
\mathcal{F} = \mathcal{F}^{t,+} \mbox{ or } \mathcal{F} =
\mathcal{F}^{t,-}. \label{label4}
\end{equation}

\begin{eg} \label{eg1} \emph{Let $F_1, ..., F_n$ be $n \geq 2$
disjoint sets, each of size at least $t$, and let $F_{n+1} =
\bigcup_{i=1}^n F_i$. Let $\mathcal{F} = \{F_1, ..., F_n,
F_{n+1}\}$. Then $\mathcal{F}^{t,+} = \{F_{n+1}\}$,
$\mathcal{F}^{t,-} = \{F_1, ..., F_n\}$ and $l(\mathcal{F},t) =
2$. By (\ref{label4}), $\beta(\mathcal{F},t) \neq
\frac{l(\mathcal{F},t)}{|\mathcal{F}|}$ (so $\beta(\mathcal{F},t)
< \frac{l(\mathcal{F},t)}{|\mathcal{F}|}$ by (\ref{label2})); one
can easily check that in fact $\beta(\mathcal{F},t) =
\beta(\mathcal{F},t,\mathcal{F}) = \frac{1}{n}$, and hence $\frac{1}{|\mathcal{F}|} < \beta(\mathcal{F},t) < \frac{l(\mathcal{F},t)}{|\mathcal{F}|}$.}
\end{eg}

Clearly, if $\mathcal{F} = \mathcal{F}^{t,+}$, then $\mathcal{F}$
is $t$-intersecting and hence $\beta(\mathcal{F},t) =
\frac{l(\mathcal{F},t)}{|\mathcal{F}|} = 1$. However, if
$\mathcal{F} = \mathcal{F}^{t,-}$, then we do not necessarily have
$\beta(\mathcal{F},t) = \frac{l(\mathcal{F},t)}{|\mathcal{F}|}$;
so the converse of (\ref{label4}) is not true.

\begin{eg} \label{eg2} \emph{Let $2 \leq m < n$, and let $F_1, ...,
F_n, F_{n+1}$ be as in Example~\ref{eg1}. Let $F_{n+2}, ..., F_{n+m}$
be sets that are disjoint from each other and from $F_{n+1} =
\bigcup_{i=1}^{n+1} F_i$. Let $\mathcal{F} = \{F_1, ...,
F_{n+m}\}$. Then $\mathcal{F} = \mathcal{F}^{t,-}$ and
$l(\mathcal{F},t) = 2$. Let $\mathcal{A} = \{F_1, ...,
F_{n+1}\}$. Then $\mathcal{A}^{t,+} = \{F_{n+1}\}$,
$\mathcal{A}^{t,-} = \{F_1, ..., F_n\}$, and
$\beta(\mathcal{F},t) \leq \beta(\mathcal{F},t,\mathcal{A}) =
\frac{1}{n} < \frac{2}{n+m} =
\frac{l(\mathcal{F},t)}{|\mathcal{F}|}$.}
\end{eg}

\subsection{The value of $\beta(\mathcal{F},t)$ for various
important families $\mathcal{F}$}\label{betasection2}
There are many important families $\mathcal{F}$ which attain the
upper bound in (\ref{label2}).

In each of the papers \cite{Borg4, Borg3, Borg2, BL2}, a
particular important family $\mathcal{F}$ is considered, and it
is proved that $|\mathcal{A}^{1,+}| +
\frac{l(\mathcal{F},1)}{|\mathcal{F}|}|\mathcal{A}^{1,-}| \leq
l(\mathcal{F},1)$ for any $\mathcal{A} \subseteq \mathcal{F}$, meaning that $\beta(\mathcal{F},1) = \frac{l(\mathcal{F},1)}{|\mathcal{F}|}$
by Proposition~\ref{beta1}. In \cite{Borg4} (a paper inspired by
\cite{H}), $\mathcal{F}$ is ${[n] \choose r}$, $r \leq n/2$ (if
$n/2< r \leq n$, then $\mathcal{F}$ is $1$-intersecting, and hence
$\beta(\mathcal{F},1) = \frac{l(\mathcal{F},1)}{|\mathcal{F}|}$
still holds). In \cite{Borg3}, $\mathcal{F}$ is
\[\mathcal{P}_{r,n} = \{\{(1,y_1), (2,y_2), ...,
(r,y_r)\} \colon y_1, y_2, ..., y_r \mbox{ are distinct elements
of } [n]\} \quad (r \in [n]),\]
which describes permutations of $r$-subsets of $[n]$ (see
\cite{Borg3}). In \cite{Borg2}, $\mathcal{F}$ is
\begin{align} \mathcal{P}_n^{(r)} = \{\{(x_1,y_1), ...,
(x_r,y_r)\} \colon &x_1, ..., x_r \mbox{ are distinct elements of
} [n], \nonumber \\
& y_1, ..., y_r \mbox{ are distinct elements of } [n]\} \quad (r \in [n]), \nonumber
\end{align}
which describes $r$-partial permutations of $[n]$ (see
\cite{Borg2}). In \cite{BL2}, $\mathcal{F}$ is the family
\[\mathcal{S}_{n,r,m} = \left\{\{(x_1,y_1), ..., (x_r,y_r)\}
\colon x_1, ..., x_r \mbox{ are distinct elements
of } [n], \, y_1, ..., y_r \in [m]\right\}\]
of $m$-signed $r$-subsets of $[n]$, where $r \in [n]$ and $m \geq 2$. For each of these families, the value of $l(\mathcal{F},1)$ is known and is attained by the intersecting sub-family $\{F \in \mathcal{F} \colon x \in F\}$ for any $x \in \bigcup_{F \in \mathcal{F}} F$ (for example, the sub-family $\left\{A \in {[n] \choose r} \colon 1 \in A \right\}$ of ${[n] \choose r}$, the sub-family $\{A \in \mathcal{P}_{r,n} \colon (1,1) \in A\}$ of $\mathcal{P}_{r,n}$, etc.); see \cite{Borg7}.

We now prove that the same holds for the \emph{power set} of a
set $X$, i.e.~the family of all subsets of $X$, which is perhaps
the most natural family one can think of. Let $2^{X}$ denote the
power set of $X$. One of the basic results in extremal set
theory is that $l(2^{[n]},1) = 2^{n-1}$ (see \cite{EKR}), and
this is generalised by our next result.

\begin{theorem} \label{betapowerset} If $\mathcal{F} = 2^{[n]}$, then
\[\beta(\mathcal{F},1) = \frac{l(\mathcal{F},1)}{|\mathcal{F}|} =
\frac{1}{2}.\]
\end{theorem}
\textbf{Proof.} Let $\mathcal{A} \subseteq \mathcal{F} = 2^{[n]}$.
Let $\mathcal{B} = \{[n] \backslash A \colon A \in
\mathcal{A}^{1,+}\}$. So $|\mathcal{B}| = |\mathcal{A}^{1,+}|$.
For any $B \in \mathcal{B}$, we have $B = [n] \backslash A$ for
some $A \in \mathcal{A}^{1,+}$, and hence, by definition of
$\mathcal{A}^{1,+}$, $B \notin \mathcal{A}$ since $|A \cap B| =
0$. So $\mathcal{A}$ and $\mathcal{B}$ are disjoint sub-families
of $\mathcal{F}$. Therefore,
\begin{equation} 2|\mathcal{A}^{1,+}| + |\mathcal{A}^{1,-}| =
|\mathcal{A}^{1,+}| + |\mathcal{B}| + |\mathcal{A}^{1,-}| =
|\mathcal{A}| + |\mathcal{B}| = |\mathcal{A} \cup \mathcal{B}|
\leq |\mathcal{F}| = 2^n \nonumber
\end{equation}
and hence, dividing throughout by $2$, we get $|\mathcal{A}^{1,+}|
+ \frac{1}{2} |\mathcal{A}^{1,-}| \leq 2^{n-1}$. It follows that a
$1$-intersecting sub-family of $\mathcal{F}$ has size at most
$2^{n-1}$ (as $\mathcal{A} = \mathcal{A}^{1,+}$ if $\mathcal{A}$
is $1$-intersecting), and this bound is attained by the trivial
$1$-intersecting sub-family $\{A \in \mathcal{F} \colon 1 \in
A\}$; so $l(\mathcal{F},1) = 2^{n-1}$ and
$\frac{l(\mathcal{F},1)}{|\mathcal{F}|} = \frac{1}{2}$. So we
have $|\mathcal{A}^{1,+}| +
\frac{l(\mathcal{F},1)}{|\mathcal{F}|} |\mathcal{A}^{1,-}| \leq
l(\mathcal{F},1)$. By Proposition~\ref{beta1}, $\beta(\mathcal{F},1) =
\frac{l(\mathcal{F},1)}{|\mathcal{F}|} =
\frac{1}{2}$.~\hfill{$\Box$}
\\

Note that by Theorems~\ref{betapowerset} and Theorem~\ref{main},
for $\mathcal{F} = 2^{[n]}$ the configuration $\mathcal{A}_1 =
... = \mathcal{A}_k = \mathcal{L}$ is optimal for both the sum
and the product for any $k \geq 2$. More precisely, we have the
following.

\begin{theorem} \label{psthm} Let $k \geq 2$, and let $\mathcal{A}_1, ..., \mathcal{A}_k$ be cross-$1$-intersecting sub-families of $2^{[n]}$. Then
\[\sum_{i = 1}^k |\mathcal{A}_i| \leq k2^{n-1} \quad
\mbox{ and } \quad \prod_{i = 1}^k |\mathcal{A}_i| \leq
2^{k(n-1)},\]
and both bounds are attained if $\mathcal{A}_1 = ... =
\mathcal{A}_k = \{A \in 2^{[n]} \colon 1 \in A\}$. Moreover, if $k > 2$,
then in both inequalities, equality holds if and only if $\mathcal{A}_1 = ... = \mathcal{A}_k = \mathcal{L}$ for some largest $1$-intersecting
sub-family $\mathcal{L}$ of $\mathcal{F}$.\footnote{At the time
of writing this paper, this result was generalised in \cite{Borg5}
for any union of power sets of sets which have a common element.}
\end{theorem}

The above results for $\beta(\mathcal{F},1)$ generalise for
$\beta(\mathcal{F},t)$ as follows. Recently, Wang and Zhang
\cite{WZ} observed that the method employed in \cite{Borg4,
Borg3, Borg2, BL2} together with a result for vertex-transitive
graphs found in \cite{AC} and also in \cite{CK} (see \cite{WZ})
immediately give us $\beta(\mathcal{F},t) =
\frac{l(\mathcal{F},t)}{|\mathcal{F}|}$ for the following very
important class of families.

We shall call a family $\mathcal{F}$ \emph{$t$-symmetric} if
there exists a group $\Gamma$ of bijections with domain
$\mathcal{F}$ and co-domain $\mathcal{F}$ such that $\Gamma$ acts
transitively on $\mathcal{F}$ and preserves the $t$-intersection
property, i.e.~for any $A, B \in \mathcal{F}$, the following hold:\\
(a) there exists $\delta \in \Gamma$ such that $B = \delta(A)$;\\
(b) if $A$ $t$-intersects $B$, then for all $\gamma \in \Gamma$, $\gamma(A)$ $t$-intersects $\gamma(B)$ .

The result proved by Wang and Zhang \cite[Corollary 2.4]{WZ} gives us the following.

\begin{theorem}[\cite{WZ}]\label{WZ} If $\mathcal{A}$ is a sub-family of a $t$-symmetric family $\mathcal{F}$, then
\begin{equation} |\mathcal{A}^{t,+}| +
\frac{l(\mathcal{F},t)}{|\mathcal{F}|}|\mathcal{A}^{t,-}| \leq
l(\mathcal{F},t). \nonumber
\end{equation}
\end{theorem}
Together with Proposition~\ref{beta1}, this immediately gives us the next result.

\begin{cor}\label{WZcor} For any $t$-symmetric family $\mathcal{F}$,
\begin{equation} \beta(\mathcal{F},t) =
\frac{l(\mathcal{F},t)}{|\mathcal{F}|}. \nonumber
\end{equation}
\end{cor}

It turns out that the families ${[n] \choose r}$,
$\mathcal{P}_{r,n}$, $\mathcal{P}_n^{(r)}$ and
$\mathcal{S}_{n,r,m}$ are $t$-symmetric for any $t \geq 1$. Also, the value of $l(\mathcal{F},t)$ is known precisely for the following cases: $\mathcal{F} = {[n] \choose r}$ (see \cite{AK1}), $\mathcal{F} = \mathcal{S}_{n,n,m}$ (see \cite{AK2, FT2}), $\mathcal{F} = \mathcal{S}_{n,r,m}$ with $n \geq (r-t+m)(t+1)/m$ (see \cite{Bey1}), $\mathcal{F} = \mathcal{P}_{n,n}$ with $n$ sufficiently larger than $t$ (see \cite{EFP}), $\mathcal{F} = \mathcal{P}_{r,n}$ with $n$ sufficiently larger than $r$ (see \cite{Borg6}), and $\mathcal{P}_n^{(r)}$ with $n$ sufficiently larger than $r$ (see \cite{Ku_thesis, Borg6, Borg}). Thus, by Corollary~\ref{WZcor}, we know $\beta(\mathcal{F},t)$ for each of these cases.

Another important family for which we have similar results is the family $\mathcal{V}_{n,r}(q)$ of all $r$-dimensional subspaces of an $n$-dimensional vector space over a $q$-element field; however, for this family, $t$-intersection is defined slightly differently, and we will discuss this separately in Section~\ref{sumproblem2}.

Now $l(2^{[n]},t)$ was determined in \cite{Kat}, and although $2^{[n]}$ is not $t$-symmetric, we will now determine $\beta(2^{[n]},t)$ using the fact that $\mathcal{S}_{n,n,2}$ is $t$-symmetric and that we also know $l(\mathcal{S}_{n,n,2},t)$ (see \cite{Kl}, and see \cite{AK2, FT2} for $\mathcal{S}_{n,n,m}$), which is in fact equal to $l(2^{[n]},t)$. Define
\[ \mathcal{K}_{n,t} =
\begin{cases}
\displaystyle \{A \subseteq [n] \colon |A| \geq (n+t)/2\} & \text{if $n+t$ is even;}\\[4mm]
\displaystyle \{A \subseteq [n] \colon |A \cap [n-1]| \geq (n+t-1)/2\} & \text{if $n+t$ is odd.}
\end{cases} \]
Katona \cite{Kat} proved that $\mathcal{K}_{n,t}$ is a largest $t$-intersecting sub-family of $2^{[n]}$ (and uniquely so up to isomorphism if $t \geq 2$); so $l(2^{[n]},t) = |\mathcal{K}_{n,t}|$. Kleitman \cite{Kl} showed that we also have $l(\mathcal{S}_{n,n,2},t) = |\mathcal{K}_{n,t}|$.

\begin{theorem}\label{betapowerset2} If $\mathcal{F} = 2^{[n]}$ and $n \geq t$, then
\[\beta(\mathcal{F},t) = \frac{l(\mathcal{F},t)}{|\mathcal{F}|} = \frac{|\mathcal{K}_{n,t}|}{2^n}.\]
\end{theorem}
\textbf{Proof.} Let $\mathcal{A} \subseteq 2^{[n]}$. For each $A \in \mathcal{A}$, let $B_A$ be the set $\{(a,1) \colon a \in A\} \cup \{(b,2) \colon b \in [n] \backslash A\}$ in $\mathcal{S}_{n,n,2}$. Let $\mathcal{B}$ be the sub-family $\{B_A \colon A \in \mathcal{A}\}$ of $\mathcal{S}_{n,n,2}$. By Theorem~\ref{WZ},
$|\mathcal{B}^{t,+}| +
\frac{l(\mathcal{S}_{n,n,2},t)}{|\mathcal{S}_{n,n,2}|}|\mathcal{B}^{t,-}| \leq
l(\mathcal{S}_{n,n,2},t)$.
Since $|\mathcal{S}_{n,n,2}| = 2^n$ and $l(\mathcal{S}_{n,n,2},t) = |\mathcal{K}_{n,t}|$, we get $|\mathcal{B}^{t,+}| +
\frac{|\mathcal{K}_{n,t}|}{2^n}|\mathcal{B}^{t,-}| \leq
|\mathcal{K}_{n,t}|$. Now we clearly have that if $A \in \mathcal{A}^{t,+}$, then $B_A \in \mathcal{B}^{t,+}$. So $|\mathcal{A}^{t,+}| = |\mathcal{B}^{t,+}| - p$ for some non-negative integer $p$, and hence, since $|\mathcal{A}^{t,+}| + |\mathcal{A}^{t,-}| = |\mathcal{A}| = |\mathcal{B}| = |\mathcal{B}^{t,+}| + |\mathcal{B}^{t,-}|$, we have $|\mathcal{A}^{t,-}| = |\mathcal{B}^{t,-}| + p$. So we have
\[|\mathcal{A}^{t,+}| + \frac{l(2^{[n]},t)}{|2^{[n]}|}|\mathcal{A}^{t,-}| = (|\mathcal{B}^{t,+}| - p) +
\frac{|\mathcal{K}_{n,t}|}{2^n} (|\mathcal{B}^{t,-}| + p)
\leq |\mathcal{B}^{t,+}| +
\frac{|\mathcal{K}_{n,t}|}{2^n}|\mathcal{B}^{t,-}| \leq
|\mathcal{K}_{n,t}|. \]
Therefore, by Proposition~\ref{beta1}, $\beta(2^{[n]},t) = \frac{l(2^{[n]},t)}{|2^{[n]}|}$ and hence the result.~\hfill{$\Box$}

\section{The maximum sum}
\label{sumproblem}

In this section we restrict our attention to the the problem of
maximising the sum of sizes of \emph{any} number of
cross-$t$-intersecting sub-families of a given family
$\mathcal{F}$. Similarly to Section~\ref{betasection}, we first
prove general results and then we provide complete solutions for
various important families.

\subsection{General results and observations}\label{sumproblem1}

We start by proving Theorem~\ref{summax2}.\\
\\
\textbf{Proof of Theorem~\ref{summax2}.} (i) is given by
Theorem~\ref{main}. Suppose $k < \kappa(\mathcal{F},t)$. So
$\frac{1}{k} > \beta(\mathcal{F},t)$.

\textit{Case 1: $\beta(\mathcal{F},t) =
\frac{l(\mathcal{F},t)}{|\mathcal{F}|}$}. So $k <
\frac{|\mathcal{F}|}{l(\mathcal{F},t)}$ and hence
$k(l(\mathcal{F},t)) < |\mathcal{F}|$. Let $\mathcal{B}_1 =
\mathcal{F}$ and $\mathcal{B}_2 = ... = \mathcal{B}_k =
\emptyset$. Since $\mathcal{B}_1, ..., \mathcal{B}_k$ are
cross-$t$-intersecting, $\sum_{i=1}^k |\mathcal{A}_i| \geq
\sum_{i=1}^k |\mathcal{B}_i|$. So we have $\sum_{i=1}^k
|\mathcal{A}_i| \geq |\mathcal{F}| > k(l(\mathcal{F},t))$.

\textit{Case 2: $\beta(\mathcal{F},t) \neq
\frac{l(\mathcal{F},t)}{|\mathcal{F}|}$}. By (\ref{label2}),
$\beta(\mathcal{F},t) < \frac{l(\mathcal{F},t)}{|\mathcal{F}|}$.
Thus, taking ${\mathcal{A}_0} \subseteq \mathcal{F}$ such that
$\beta(\mathcal{F},t,\mathcal{A}_0) = \beta(\mathcal{F},t)$, we
have ${\mathcal{A}_0}^{t,-} \neq \emptyset$ and
$|{\mathcal{A}_0}^{t,+}| +
\beta(\mathcal{F},t)|{\mathcal{A}_0}^{t,-}| = l(\mathcal{F},t)$.
Since $\frac{1}{k} > \beta(\mathcal{F},t)$,
$|{\mathcal{A}_0}^{t,+}| + \frac{1}{k}|{\mathcal{A}_0}^{t,-}| >
l(\mathcal{F},t)$. Let $\mathcal{B}_1 = \mathcal{A}_0$ and
$\mathcal{B}_2 = ... = \mathcal{B}_k = {\mathcal{A}_0}^{t,+}$.
Then
\[\sum_{i=1}^k |\mathcal{B}_i| = (|{\mathcal{A}_0}^{t,-}| +
|{\mathcal{A}_0}^{t,+}|) + (k-1)|{\mathcal{A}_0}^{t,+}| = k\left(
|{\mathcal{A}_0}^{t,+}| +
\frac{1}{k}|{\mathcal{A}_0}^{t,-}|\right)\]
and hence $\sum_{i=1}^k |\mathcal{B}_i| > k(l(\mathcal{F},t))$.
Since $\mathcal{B}_1, ..., \mathcal{B}_k$ are
cross-$t$-intersecting, we have $\sum_{i=1}^k |\mathcal{A}_i| \geq
\sum_{i=1}^k |\mathcal{B}_i| >
k(l(\mathcal{F},t))$.~\hfill{$\Box$}\\

As we have seen in Section~\ref{betasection}, the case
$\beta(\mathcal{F},t) = \frac{l(\mathcal{F},t)}{|\mathcal{F}|}$
deserves very special attention. For this particularly interesting
case, we have the following precise result, which gives us the
maximum sum of sizes for any $k \geq 2$, and characterises optimal configurations. Recall that $\max\{|F| \colon F \in \mathcal{F}\}$ is denoted by $\alpha(\mathcal{F})$.

\begin{theorem}\label{summax3} Let $\mathcal{F}$ be a family with
$\alpha(\mathcal{F}) \geq t$ and $\beta(\mathcal{F},t) =
\frac{l(\mathcal{F},t)}{|\mathcal{F}|}$. Let $\mathcal{A}_1, ...,
\mathcal{A}_k$ be cross-$t$-intersecting sub-families of
$\mathcal{F}$ such that $\sum_{i=1}^k |\mathcal{A}_i|$ is maximum. Then
\[ \sum_{i=1}^k |\mathcal{A}_i| =
\begin{cases}
\displaystyle |\mathcal{F}| & \text{if $k \leq
\frac{|\mathcal{F}|}{l(\mathcal{F},t)}$;}\\[4mm]
\displaystyle k(l(\mathcal{F},t)) & \text{if $k \geq
\frac{|\mathcal{F}|}{l(\mathcal{F},t)}$.}
\end{cases} \]
Moreover,  \\
(i) if $k < \frac{|\mathcal{F}|}{l(\mathcal{F},t)}$, then
$\mathcal{A}_i = {\mathcal{A}_i}^{t,-}$ for all $i \in [k]$, and
$\mathcal{A}_1, ..., \mathcal{A}_k$ partition $\mathcal{F}$; \\
(ii) if $k > \frac{|\mathcal{F}|}{l(\mathcal{F},t)}$, then
$\mathcal{A}_1 = ... = \mathcal{A}_k = \mathcal{L}$ for some
largest $t$-intersecting sub-family $\mathcal{L}$ of
$\mathcal{F}$.
\end{theorem}

\begin{rem} \label{summaxrem}\emph{An optimal configuration
for the case $k \leq \frac{|\mathcal{F}|}{l(\mathcal{F},t)}$ is
the one with $\mathcal{A}_1 = \mathcal{F}$ and $\mathcal{A}_2 =
... = \mathcal{A}_k = \emptyset$; we will call this the
\emph{trivial configuration}. If $k =
\frac{|\mathcal{F}|}{l(\mathcal{F},t)}$, then the configuration
$\mathcal{A}_1 = ... = \mathcal{A}_k = \mathcal{L}$ is also
optimal. For each of the cases $k <
\frac{|\mathcal{F}|}{l(\mathcal{F},t)}$ and $k =
\frac{|\mathcal{F}|}{l(\mathcal{F},t)}$, it is possible to have
other optimal configurations but it is also possible to not have
any others; \cite[Theorem~1.4]{Borg3} gives an example of each
of these possibilities for $t=1$.}
\end{rem}
\textbf{Proof of Theorem~\ref{summax3}}. Since
$\beta(\mathcal{F},t) = \frac{l(\mathcal{F},t)}{|\mathcal{F}|}$,
$\kappa(\mathcal{F},t) = \frac{|\mathcal{F}|}{l(\mathcal{F},t)}$.
So the case $k \geq \frac{|\mathcal{F}|}{l(\mathcal{F},t)}$ is
given by Theorem~\ref{main}.

Let $\mathcal{A} = \bigcup_{i = 1}^k \mathcal{A}_i$.
Lemma~\ref{intxint} tells us that $\mathcal{A}^{t,+} =
\bigcup_{i=1}^k {\mathcal{A}_i}^{t,+}$, $\mathcal{A}^{t,-} =
\bigcup_{i = 1}^k {\mathcal{A}_i}^{t,-}$, and
${\mathcal{A}_1}^{t,-}, ..., {\mathcal{A}_k}^{t,-}$ partition
${\mathcal{A}}^{t,-}$. So we have
\begin{align} \sum_{i=1}^k |\mathcal{A}_i| &= \sum_{i=1}^k
|{\mathcal{A}_i}^{t,-}| + \sum_{i=1}^k |{\mathcal{A}_i}^{t,+}|
\leq |\mathcal{A}^{t,-}| + k|\mathcal{A}^{t,+}| \nonumber \\
&\leq \frac{1}{\beta(\mathcal{F},t)}(l(\mathcal{F},t) -
|\mathcal{A}^{t,+}|) + k|\mathcal{A}^{t,+}| \quad \quad \mbox{(by
(\ref{label1}))} \nonumber \\
&= \frac{|\mathcal{F}|}{l(\mathcal{F},t)}(l(\mathcal{F},t) -
|\mathcal{A}^{t,+}|) + k|\mathcal{A}^{t,+}| = |\mathcal{F}| +
\left(k - \frac{|\mathcal{F}|}{l(\mathcal{F},t)} \right)
|\mathcal{A}^{t,+}| \label{main3}
\end{align}
(note that the case $k \geq
\frac{|\mathcal{F}|}{l(\mathcal{F},t)}$ can be deduced from
(\ref{main3}) since $|\mathcal{A}^{t,+}| \leq l(\mathcal{F},t)$).
Suppose $k \leq \frac{|\mathcal{F}|}{l(\mathcal{F},t)}$. Then
$\sum_{i=1}^k |\mathcal{A}_i| \leq |\mathcal{F}|$, and if $k <
\frac{|\mathcal{F}|}{l(\mathcal{F},t)}$, then, by (\ref{main3}),
the bound is attained only if $\mathcal{A}^{t,+} = \emptyset$ and
$\mathcal{A}^{t,-} = \mathcal{F}$, which implies that
$\mathcal{A}_i = {\mathcal{A}_i}^{t,-}$ for all $i \in [k]$, and that
$\mathcal{A}_1, ..., \mathcal{A}_k$ partition $\mathcal{F}$. Now
let $\mathcal{B}_1 = \mathcal{F}$ and $\mathcal{B}_2 = ... =
\mathcal{B}_k = \emptyset$. Since $\mathcal{B}_1, ...,
\mathcal{B}_k$ are cross-$t$-intersecting and $\sum_{i=1}^k
|\mathcal{B}_i| = |\mathcal{F}|$, we have $|\mathcal{F}| \leq
\sum_{i=1}^k |\mathcal{A}_i|$ (as $\sum_{i=1}^k |\mathcal{A}_i|$
is maximum). Together with $\sum_{i=1}^k |\mathcal{A}_i| \leq
|\mathcal{F}|$, this gives us $\sum_{i=1}^k |\mathcal{A}_i| =
|\mathcal{F}|$. Hence the result.~\hfill{$\Box$} \\

The results above raise the following question: can we say
something in general about the structure of an optimal
configuration for $k < \kappa(\mathcal{F},t)$? An answer is given
by the next result, which in particular describes an optimal
configuration.

\begin{prop} \label{summax} Let $\mathcal{F}$ and
$\mathcal{A}_1, ..., \mathcal{A}_k$ be as in
Theorem~\ref{summax2}. Let $\mathcal{A} = \bigcup_{i=1}^k
\mathcal{A}_i$. Let $\mathcal{A}_0$ be a sub-family of
$\mathcal{F}$ such that $|{\mathcal{A}_0}^{t,+}| +
\frac{1}{k}|{\mathcal{A}_0}^{t,-}|$ is maximum, and let
$\mathcal{B}_1 = \mathcal{A}_0$ and $\mathcal{B}_2 = ... =
\mathcal{B}_k = {\mathcal{A}_0}^{t,+}$. Then \\
(i) $\mathcal{B}_1, ..., \mathcal{B}_k$ are
cross-$t$-intersecting sub-families of $\mathcal{F}$,\\
(ii) $\sum_{i=1}^k |\mathcal{A}_i| = \sum_{i=1}^k
|\mathcal{B}_i|$ and $|{\mathcal{A}}^{t,+}| +
\frac{1}{k}|{\mathcal{A}}^{t,-}| = |{\mathcal{A}_0}^{t,+}| +
\frac{1}{k}|{\mathcal{A}_0}^{t,-}|$.
\end{prop}
\textbf{Proof.} (i) is trivial. As in the proof of
Theorem~\ref{main}, $\sum_{i=1}^k |\mathcal{A}_i| \leq
k\left(|\mathcal{A}^{t,+}| +
\frac{1}{k}|\mathcal{A}^{t,-}|\right)$. Thus, by the choice of
$\mathcal{A}_0$, $\sum_{i=1}^k |\mathcal{A}_i| \leq
k\left(|{\mathcal{A}_0}^{t,+}| +
\frac{1}{k}|{\mathcal{A}_0}^{t,-}| \right) = \sum_{i=1}^k
|\mathcal{B}_i|$, where the equality follows as in the proof of
Theorem~\ref{summax2}. Now by (i) and the choice of
$\mathcal{A}_1, ..., \mathcal{A}_k$, $\sum_{i=1}^k
|\mathcal{A}_i| \geq \sum_{i=1}^k |\mathcal{B}_i|$. So we
actually have $\sum_{i=1}^k |\mathcal{A}_i| = \sum_{i=1}^k
|\mathcal{B}_i|$ and hence $|{\mathcal{A}}^{t,+}| +
\frac{1}{k}|{\mathcal{A}}^{t,-}| = |{\mathcal{A}_0}^{t,+}| +
\frac{1}{k}|{\mathcal{A}_0}^{t,-}|$.~\hfill{$\Box$}

\begin{rem} \label{rem2} \emph{We know from Theorems~\ref{main}
and \ref{summax2} that the configuration $\mathcal{A}_1, ...,
\mathcal{A}_k = \mathcal{L}$ is always optimal for $k \geq
\kappa(\mathcal{F},t)$ (and uniquely so if $k >
\kappa(\mathcal{F},t)$) and never optimal for $k <
\kappa(\mathcal{F},t)$. Theorem~\ref{summax3} tells us that the
trivial configuration (see Remark~\ref{summaxrem}) is always
optimal for $k < \kappa(\mathcal{F},t)$ if $\beta(\mathcal{F},t) =
\frac{l(\mathcal{F},t)}{|\mathcal{F}|}$. However, as
Proposition~\ref{summax} suggests, the trivial configuration may
not be optimal for $k < \kappa(\mathcal{F},t)$ if
$\beta(\mathcal{F},t) \neq
\frac{l(\mathcal{F},t)}{|\mathcal{F}|}$, meaning that it is
possible to have $\mathcal{F}$ and $k$ for which neither of the
two simple configurations mentioned in Remark~\ref{summaxrem}
give a maximum sum. }
\end{rem}

\begin{prop} \label{summax4} Let $\mathcal{F}$ be a family with
$\alpha(\mathcal{F}) \geq t$, $\mathcal{F}^{t,+} \neq \emptyset$
and $\mathcal{F}^{t,-} \neq \emptyset$. Suppose $2 \leq k <
\kappa(\mathcal{F},t)$ and $\mathcal{A}_1, ..., \mathcal{A}_k$
are cross-$t$-intersecting sub-families of $\mathcal{F}$ such
that $\sum_{i=1}^k |\mathcal{A}_i|$ is maximum. Then we neither
have $\mathcal{A}_i = \mathcal{F}$ for some $i \in [k]$ and
$\mathcal{A}_j = \emptyset$ for all $j \in [k] \backslash \{i\}$
nor have $\mathcal{A}_1 = ... = \mathcal{A}_k = \mathcal{L}$ for
some largest $t$-intersecting sub-family $\mathcal{L}$ of
$\mathcal{F}$.
\end{prop}
\textbf{Proof.} By Theorem~\ref{summax2}, we do not have
$\mathcal{A}_1 = ... = \mathcal{A}_k = \mathcal{L}$. Let $\mathcal{B}_1 = \mathcal{F}$ and
$\mathcal{B}_2 = ... = \mathcal{B}_k = \mathcal{F}^{t,+}$. Since
$\mathcal{B}_1, ..., \mathcal{B}_k$ are cross-$t$-intersecting,
$\sum_{i=1}^k |\mathcal{A}_i| \geq \sum_{i=1}^k |\mathcal{B}_i|$.
So $\sum_{i=1}^k |\mathcal{A}_i| \geq |\mathcal{F}| +
(k-1)|\mathcal{F}^{t,+}| > |\mathcal{F}|$. The result
follows.~\hfill{$\Box$} \\

Note that if $\mathcal{F}$ is as in the above proposition, then for all $k \geq 2$, the trivial configuration is not optimal; this is immediate from the proof of the proposition.

An example of a family $\mathcal{F}$ as in the above result is the
one in Example~\ref{eg1}. The example below shows that the
phenomenon described at the end of Remark~\ref{rem2} may also
happen when $\mathcal{F}^{t,+} = \emptyset$ and hence
$\mathcal{F} = \mathcal{F}^{t,-}$ (it cannot happen when
$\mathcal{F}^{t,-} = \emptyset$, because then $\mathcal{F}$
itself is $t$-intersecting and hence we can take $\mathcal{A}_1 =
... = \mathcal{A}_k = \mathcal{F}$).

\begin{eg} \label{eg3} \emph{Let $2 \leq m < k < n$. Let
$\mathcal{F} = \{F_1, ..., F_{n+m}\}$ be as in Example~\ref{eg2}.
Let $\mathcal{A}_1 = \{F_1, ...., F_{n+1}\}$ and $\mathcal{A}_2 =
...= \mathcal{A}_k = \{F_{n+1}\}$. Then $\mathcal{A}_1, ...,
\mathcal{A}_k$ are cross-$t$-intersecting and $\sum_{i=1}^k
|\mathcal{A}_i| = n+k > \max\{n+m, 2k\} = \max\{|\mathcal{F}|,
k(l(\mathcal{F},t))\}$.}
\end{eg}
However, unlike Proposition~\ref{summax4}, if $\mathcal{F} =
\mathcal{F}^{t,-}$, $\beta(\mathcal{F},t) \neq
\frac{l(\mathcal{F},t)}{|\mathcal{F}|}$ and $k <
\kappa(\mathcal{F},t)$, then the trivial configuration may still
give a maximum sum.

\begin{eg} \label{eg4} \emph{Let $\mathcal{F} =
\{F_1, ..., F_{n+m}\}$ be as in Example~\ref{eg2}, and let $2 \leq
k \leq m$. We have $\mathcal{F} = \mathcal{F}^{t,-}$. If
$\mathcal{A} \subseteq \mathcal{F}$ and $\mathcal{A}^{t,+} \neq
\emptyset$, then one of the following holds: (i) $|\mathcal{A}| =
1$, (ii) $\mathcal{A} = \mathcal{A}^{t,+} = \{F_i, F_{n+1}\}$ for
some $i \in [n]$, (iii) $\mathcal{A}^{t,+} = \{F_{n+1}\}$ and
$\mathcal{A}^{t,-} \subseteq \{F_1, ..., F_n\}$. It is therefore
easy to check that $\beta(\mathcal{F},t) = \beta(\mathcal{F},t,
\{F_1, ..., F_{n+1}\}) = \frac{1}{n} < \frac{2}{n+m} =
\frac{l(\mathcal{F},t)}{|\mathcal{F}|}$. Since $k \leq m < n$, $k
< \kappa(\mathcal{F},t)$. It is also easy to check that if
$\mathcal{A} \subseteq \mathcal{F}$, then $|\mathcal{A}^{t,+}| +
\frac{1}{k}|\mathcal{A}^{t,-}|$ is maximum if $\mathcal{A} =
\mathcal{F}$. By Proposition~\ref{summax}, the trivial
configuration gives a maximum sum.}
\end{eg}

\subsection{Solutions for various important
families}\label{sumproblem2}

Section~\ref{betasection2} gives the values of
$\beta(\mathcal{F},t)$ that we know for the families $2^{[n]}$,
${[n] \choose r}$, $\mathcal{P}_{r,n}$, $\mathcal{P}_n^{(r)}$ and
$\mathcal{S}_{n,r,m}$, and they all turn out to be the maximum
possible value $\frac{l(\mathcal{F},t)}{|\mathcal{F}|}$. Thus, by
Theorem~\ref{summax3}, for all these cases we know the maximum
sum of sizes of any $k \geq 2$ cross-$t$-intersecting
sub-families of $\mathcal{F}$, and we also know that at least one
of the trivial configuration (see Remark~\ref{summaxrem}) and
the configuration $\mathcal{A}_1 = ... = \mathcal{A}_k =
\mathcal{L}$ is optimal, with the latter being the unique optimal
configuration when $k > \frac{|\mathcal{F}|}{l(\mathcal{F},t)}$.
We point out that \cite[Theorem 2.5]{WZ} tells us that in addition
to this, for the cases we are discussing except for the one with $\mathcal{F} = \mathcal{P}_{3,3}$ and $t=1$, when $k < \frac{|\mathcal{F}|}{l(\mathcal{F},t)}$ the trivial configuration
is the unique optimal configuration if we simply insist that
$\mathcal{A}_1 \neq \emptyset$ (for $2^{[n]}$, this emerges from the correspondence with $\mathcal{S}_{n,n,2}$ used in the proof of Theorem~\ref{betapowerset2}). As pointed out in
Remark~\ref{summaxrem}, for the case $k =
\frac{|\mathcal{F}|}{l(\mathcal{F},t)}$ there may be other other
optimal configurations apart from the two mentioned above.

Now recall the family $\mathcal{V}_{n,r}(q)$ defined in Secion~\ref{betasection2}. If $A, B \in \mathcal{V}_{n,r}(q)$ such that $\dim(A \cap B) \geq t$, then, with slight abuse of terminology, we say that $A$ \emph{$t$-intersects} $B$. For $\mathcal{V}_{n,r}(q)$, we work with this definition of $t$-intersection instead of the usual one, and so we define $t$-intersecting sub-families, cross-$t$-intersecting sub-families, $l(\mathcal{V}_{n,r}(q),t)$, $\beta(\mathcal{V}_{n,r}(q),t)$, and so on, accordingly. The value of $l(\mathcal{V}_{n,r}(q),t)$ was determined in \cite{FW}. $\mathcal{V}_{n,r}(q)$ is $t$-symmetric; see \cite[Example~1.3]{WZ}. By \cite[Corollary~2.4]{WZ}, the statement of Theorem~\ref{WZ} holds for $\mathcal{V}_{n,r}(q)$. Thus, by applying the argument in the proof of Proposition~\ref{beta1} to $\mathcal{V}_{n,r}(q)$, we obtain $\beta(\mathcal{V}_{n,r}(q),t) = \frac{l(\mathcal{V}_{n,r}(q),t)}{|\mathcal{V}_{n,r}(q)|}$. \cite[Theorem 2.5]{WZ} solved the problem of maximising the sum of sizes of $k \geq 2$ cross-$t$-intersecting sub-families of $\mathcal{V}_{n,r}(q)$.

\section{The maximum product}
\label{productproblem}

In this section we restrict our attention to the the problem of
maximising the product of sizes of cross-$t$-intersecting
sub-families of a given family $\mathcal{F}$. Similarly to the
two preceding sections, we first reveal some interesting facts and
then we provide solutions for various important families.

\subsection{General results and observations}
\label{productproblem1}

Theorem~\ref{main} tells us that the configuration $\mathcal{A}_1
= ... = \mathcal{A}_k = \mathcal{L}$ (where $\mathcal{L}$ is a
largest $t$-intersecting sub-family of $\mathcal{F}$) gives both
a maximum sum and a maximum product of sizes when $k \geq
\kappa(\mathcal{F},t)$, and Theorem~\ref{summax2} tells us that
this configuration never gives a maximum sum when $k <
\kappa(\mathcal{F},t)$. However, this configuration may still
give a maximum product when $k < \kappa(\mathcal{F},t)$. For
example, the main result in \cite{MT} tells us that the product
of sizes of $2$ cross-$1$-intersecting sub-families
$\mathcal{A}_1$ and $\mathcal{A}_2$ of ${[n] \choose r}$ is maximum if $\mathcal{A}_1 = \mathcal{A}_2 = \mathcal{L}$ for some
largest $1$-intersecting sub-family $\mathcal{L}$ of ${[n]
\choose r}$, where $\mathcal{L} = {[n] \choose r}$ if $n/2 < r \leq n$,
and by the classical result in \cite{EKR}, $\mathcal{L}$ is of
size ${n-1 \choose r-1}$ (the size of the
$1$-intersecting sub-family $\{A \in {[n] \choose r} \colon 1 \in
A\}$ of ${[n] \choose r}$) if $r \leq n/2$; note that if $n >
2r$, then, since $\beta\left( {[n] \choose r}, 1 \right) =
\frac{|\mathcal{L}|}{{n \choose r}} = \frac{r}{n}$ (see
Section~\ref{betasection2}), we have $2 < \kappa\left( {[n]
\choose r}, 1 \right)$. The general cross-$t$-intersection version (also for $2$ sub-families) is given in \cite{Tok} for $n$ sufficiently large; see Theorem~\ref{Tokcross}. Other results of this kind are given in the next sub-section. The following tells us that such product results generalise to $k$ sub-families for any $k \geq 2$.

\begin{lemma}\label{prodext} Let $\mathcal{L}$ be a largest
$t$-intersecting sub-family of a family $\mathcal{F}$. Suppose
that the product of sizes of $p$ cross-$t$-intersecting
sub-families $\mathcal{B}_1, ..., \mathcal{B}_p$ of $\mathcal{F}$
is maximum if $\mathcal{B}_1 = ... = \mathcal{B}_p =
\mathcal{L}$. Then for any $k \geq p$, the product of sizes of $k$
cross-$t$-intersecting sub-families $\mathcal{A}_1, ...,
\mathcal{A}_k$ of $\mathcal{F}$ is maximum if $\mathcal{A}_1 =
... = \mathcal{A}_k = \mathcal{L}$.
\end{lemma}
We first prove the following result, which immediately yields the above result.

\begin{lemma} \label{prodextlemma} Let $k \geq p$, and let $x_1,
..., x_k, y_1, ..., y_k$ be non-negative real numbers such that
$\prod_{i \in I} x_i \leq \prod_{i \in I} y_i$ for any subset $I$
of $[k]$ of size $p$. Then $\prod_{i=1}^k x_i \leq \prod_{i = 1}^k
y_i$.
\end{lemma}
\textbf{Proof.} Let mod$^*$ represent the usual \emph{modulo
operation} with the exception that for any two positive integers $a$ and $b$, $ba \; {\rm mod}^* \; a$ is $a$ instead of $0$. We have
\begin{equation} \left( \prod_{i=1}^k x_i \right)^p = \prod_{i = 0}^{k-1} \prod_{j = 1}^p x_{(ip+j) \; {\rm mod}^* \; k} \leq \prod_{i = 0}^{k-1} \prod_{j = 1}^p y_{(ip+j) \; {\rm mod}^* \; k} = \left( \prod_{i=1}^k y_i \right)^p. \nonumber
\end{equation}
Hence the result.~\hfill{$\Box$} \\
\\
\textbf{Proof of Lemma~\ref{prodext}.} By our assumption,
$\prod_{i \in I} |\mathcal{A}_i| \leq \left( l(\mathcal{F},t)
\right)^p$ for any subset $I$ of $[k]$ of size $p$. By
Lemma~\ref{prodextlemma} with $x_i = |\mathcal{A}_i|$ and $y_i =
l(\mathcal{F},t)$ for all $i \in [k]$, $\prod_{i = 1}^k
|\mathcal{A}_i| \leq \left(l(\mathcal{F},t) \right)^k$. The
result follows.~\hfill{$\Box$} \\

We now prove Remark~\ref{productthreshold}. More precisely, we
will show that for any $t \geq 1$ and any $p \geq 3$, there are
families $\mathcal{F}$ with $\kappa(\mathcal{F},t) = p$ such
that, unlike the case when $k \geq \kappa(\mathcal{F},t)$ (see
Theorem~\ref{main}), for any $2 \leq k < \kappa(\mathcal{F},t)$,
the product of sizes of $k$ cross-$t$-intersecting sub-families
$\mathcal{A}_1, ..., \mathcal{A}_k$ of $\mathcal{F}$ is not maximum if $\mathcal{A}_1 = ... = \mathcal{A}_k = \mathcal{L}$
for some largest $t$-intersecting sub-family $\mathcal{L}$ of
$\mathcal{F}$. Our aim is to construct a family $\mathcal{P}$ of size $p^2$ that can be partitioned into $p$ cross-$1$-intersecting families $\mathcal{P}_1, ..., \mathcal{P}_p$, each of size $p$, such that for any $i \in [p]$, the $p$ sets $A_{i,1}, ..., A_{i,p}$ in $\mathcal{P}_i$ are disjoint. Then we take $\mathcal{B}$ to be the family obtained from $\mathcal{P}$ by replacing each element $u$ of the union of all sets in $\mathcal{P}$ by $t$ new elements $u_1, ..., u_t$.

\begin{construction} \label{geomcon} \emph{Let $p \geq 3$ be
an integer. Let $m_1, ..., m_p$ and $c_1, ..., c_p$ be distinct
real numbers. For any $i, j \in [p]$, let $L_{i,j}$ be the
straight line in $\mathbb{R}^2$ arising from the function $y
\colon \mathbb{R} \rightarrow \mathbb{R}$ defined by $y(x) = m_ix
+ c_j$. For any $i, j \in [p]$, let $A_{i,j}$ be the set of
all points (i.e.~co-ordinates) of intersection of $L_{i,j}$ with
the other lines $L_{i',j'}$, i.e.
\[A_{i,j} = \left\{(a,b) \in \mathbb{R}^2 \colon
\exists \, i', j' \in [p], \, (i',j') \neq (i,j), \, \mbox{such
that } L_{i,j} \mbox{ intersects } L_{i',j'} \mbox{ at }
(a,b) \right\}.\]
Let $(a_1, b_1), ..., (a_s,b_s)$ be the distinct co-ordinates in
the set $\bigcup_{i=1}^p\bigcup_{j=1}^p A_{i,j}$ of all points of pairwise intersection of these lines, and let $T_{(a_1,b_1)}, ...,
T_{(a_s,b_s)}$ be disjoint sets of size $t$. For any $i, j \in
[p]$, let $B_{i,j} = \bigcup_{(a,b) \in A_{i,j}} T_{(a,b)}$; so
$B_{i,j}$ is simply the set obtained by replacing each point
$(a,b)$ in $A_{i,j}$ by the $t$ elements of the corresponding set $T_{(a,b)}$. For each $i \in [p]$, let $\mathcal{B}_i = \{B_{i,1}, ...,
B_{i,p}\}$. Now let $\mathcal{B} = \bigcup_{i=1}^p \mathcal{B}_i =
\{B_{i,j} \colon i, j \in [p]\}$.}
\end{construction}

\begin{theorem} \label{geomthm} Let $\mathcal{B}$ be as in
Construction~\ref{geomcon}. Let $\mathcal{L}$ be a largest
$t$-intersecting sub-family of $\mathcal{B}$, and let
$\mathcal{A}_1, ..., \mathcal{A}_k$ be cross-$t$-intersecting
sub-families of $\mathcal{B}$. Then:\\
(i) $\kappa(\mathcal{B},t) = |\mathcal{L}| = p$;\\
(ii) if $k \geq \kappa(\mathcal{B},t)$ and $\mathcal{A}_1 = ... =
\mathcal{A}_k = \mathcal{L}$, then $\prod_{i=1}^k
|\mathcal{A}_i|$ is maximum; \\
(iii) if $k < \kappa(\mathcal{B},t)$ and $\mathcal{A}_1 = ... =
\mathcal{A}_k = \mathcal{L}$, then $\prod_{i=1}^k
|\mathcal{A}_i|$ is not maximum.
\end{theorem}
\textbf{Proof.} Let $\mathcal{I}$ be a $t$-intersecting sub-family
of $\mathcal{B}$. For each $i \in [p]$, the lines $L_{i,1}, ...,
L_{i,p}$ have the same gradient $m_i$, and hence, since $c_1,
..., c_p$ are distinct, $L_{i,1}, ..., L_{i,p}$ are distinct
parallel lines, meaning that no two intersect. Thus, for each $i \in
[p]$, $\mathcal{I}$ contains at most one of the sets in
$\mathcal{B}_i$. So $|\mathcal{I}| \leq p$. Now for any $i, i',
j, j' \in [p]$ with $i \neq i'$, $L_{i,j}$ intersects $L_{i',j'}$
(at one point) since $m_i \neq m_{i'}$. So $\{A_{1,1}, A_{2,1},
..., A_{p,1}\}$ is a $1$-intersecting family (in fact, $(0,c_1) \in A_{i,1}$ for each $i \in [p]$) of size $p$, meaning
that $\{B_{1,1}, B_{2,1}, ..., B_{p,1}\}$ is a $t$-intersecting
sub-family of $\mathcal{B}$ of size $p$, and hence a largest $t$-intersecting sub-family of $\mathcal{B}$. So $|\mathcal{L}| = p = l(\mathcal{B},t)$.

Let $\mathcal{A}$ be a sub-family of $\mathcal{B}$. If
$\mathcal{A}^{t,-} = \emptyset$, then $\beta(\mathcal{B},t,\mathcal{A}) = \frac{l(\mathcal{B},t)}{|\mathcal{B}|}$. Suppose
$\mathcal{A}^{t,-} \neq \emptyset$. By the same argument for
$\mathcal{I}$ above, for each $i \in [p]$, $\mathcal{A}^{t,+}$
contains at most one of the sets in $\mathcal{B}_i$, and if it
does contain one of these sets, then, by definition of
$\mathcal{A}^{t,+}$, $\mathcal{A}$ contains no other set in
$\mathcal{B}_i$. Let $S = \{i \in [p] \colon \mathcal{A}^{t,+}
\mbox{ contains one of the sets in } \mathcal{B}_i\}$. Then
$|\mathcal{A}^{t,+}| = |S|$, and for each $s \in S$,
$\mathcal{A}^{t,-}$ contains no set in $\mathcal{B}_s$. So
$\mathcal{A}^{t,-} \subseteq \bigcup_{j \in [p] \backslash S}
\mathcal{B}_j$ and hence $|\mathcal{A}^{t,-}| \leq (p - |S|)p$.
Note that $|S| < p$ since $|\mathcal{A}^{t,-}| > 0$. So we have
\[\beta(\mathcal{B},t,\mathcal{A}) = \frac{l(\mathcal{B},t) -
|\mathcal{A}^{t,+}|}{|\mathcal{A}^{t,-}|} \geq \frac{p - |S|}{(p -
|S|)p} = \frac{1}{p} = \frac{p}{p^2} =
\frac{l(\mathcal{B},t)}{|\mathcal{B}|} =
\beta(\mathcal{B},t,\mathcal{L}).\]
Therefore, $\beta(\mathcal{B},t) = \frac{1}{p}$ and hence
$\kappa(\mathcal{B},t) = p$. Hence (i).\medskip
\\
Part (ii) is given by Theorem~\ref{main}.\medskip
\\
Finally, suppose $k < \kappa(\mathcal{B},t)$. So $p \geq k+1$. Let
$\mathcal{B}_k' = \bigcup_{i=k}^p \mathcal{B}_i$. As we mentioned
above, for $i \neq i'$, any two lines $L_{i,j}$ and $L_{i',j'}$
intersect on a point $(a,b)$, and hence the $t$-set $T_{(a,b)}$
is a subset of $B_{i,j} \cap B_{i',j'}$.  So $\mathcal{B}_1, ...,
\mathcal{B}_{k-1}, \mathcal{B}_k'$ are cross-$t$-intersecting
sub-families of $\mathcal{B}$, and the product of their sizes is
$p^{k-1}(p-k+1)p > p^k = |\mathcal{L}|^k$. Hence
(iii).~\hfill{$\Box$}

\subsection{Solutions for various important
families}\label{productproblem2}

The cross-$t$-intersection problem for the product is more
difficult than that for the sum, and hence less is known about the
product. However, various breakthroughs have been made for the
special families in Section~\ref{betasection2}.

Consider first the family $2^{[n]}$. For $t=1$ we have the
complete solution given by Theorem~\ref{psthm}, and for $t \geq 1$
we have the following.

\begin{theorem} [\cite{MT2}] Let $\mathcal{A}_1$ and $\mathcal{A}_2$ be cross-$t$-intersecting sub-families of $2^{[n]}$, where $1 \leq t \leq n$. Let $\mathcal{K}_{1} = \{A \subseteq
[n] \colon |A| \geq (n+t)/2\}$, $\mathcal{K}_{2} = \{A \subseteq
[n] \colon |A \cap [n-1]| \geq (n+t-1)/2\}$ and $\mathcal{K}_{3} = \{A \subseteq [n] \colon |A| \geq (n+t-1)/2\}$.\\
(i) If $n + t$ is even, then $|\mathcal{A}_1| |\mathcal{A}_2|\leq |\mathcal{K}_{1}|^2$. \\
(ii) If $n + t$ is odd, then $|\mathcal{A}_1|
|\mathcal{A}_2| \leq \max\{|\mathcal{K}_{2}|^2, |\mathcal{K}_{1}||\mathcal{K}_{3}|\}$.
\end{theorem}
Thus, by Lemma~\ref{prodext}, if $n+t$ is even, then the
product of $k \geq 2$ cross-$t$-intersecting sub-families
$\mathcal{A}_1, ..., \mathcal{A}_k$ of $2^{[n]}$ is maximum if
$\mathcal{A}_1 = ... = \mathcal{A}_k = \mathcal{K}_{1}$;
however, it is not known what the maximum product is when $n+t$
is odd and $k \geq 3$.

The following theorems were proved for $2$ sub-families, and for each one of them, we obtain the generalisation to any $k \geq 2$ sub-families from Lemma~\ref{prodext} (with $p = 2$).

For the family ${[n] \choose r}$, we have the next two results.

\begin{theorem}[\cite{Pyber,MT}] Let $\mathcal{A}_1$ and $\mathcal{A}_2$ be cross-$1$-intersecting sub-families of ${[n] \choose r}$, where $1 \leq r \leq n/2$. Then
\[|\mathcal{A}_1||\mathcal{A}_2| \leq {n-1 \choose r-1}^2, \]
and equality holds if $\mathcal{A}_1 = \mathcal{A}_2 = \left\{A \in {[n] \choose r} \colon 1 \in A\right\}$.
\end{theorem}
The result for $n/2 < r \leq n$ is trivial; in this case, $\mathcal{A}_1$ and $\mathcal{A}_2$ are cross-$1$-intersecting if each one of them is the whole family ${[n] \choose r}$.

\begin{theorem}[\cite{Tok}] \label{Tokcross} Let $\mathcal{A}_1$ and $\mathcal{A}_2$ be cross-$t$-intersecting sub-families of ${[n] \choose r}$, where $1 \leq t \leq r$. If $n$ is sufficiently large, then
\[|\mathcal{A}_1||\mathcal{A}_2| \leq {n-t \choose r-t}^2,\]
and equality holds if $\mathcal{A}_1 = \mathcal{A}_2 =
\left\{A \in {[n] \choose r} \colon [t] \subset A \right\}$.
\end{theorem}

Finally, for $\mathcal{P}_{n,n}$ we have the next two results.

\begin{theorem}[\cite{E}] Let $\mathcal{A}_1$ and $\mathcal{A}_2$ be
cross-$1$-intersecting sub-families of $\mathcal{P}_{n,n}$, where $n \geq 4$. Then
\[|\mathcal{A}_1||\mathcal{A}_2| \leq ((n-1)!)^2, \]
and equality holds if $\mathcal{A}_1 = \mathcal{A}_2 =
\{A \in \mathcal{P}_{n,n} \colon (1,1) \in A\}$.
\end{theorem}

\begin{theorem} [\cite{EFP}] Let $\mathcal{A}_1$ and $\mathcal{A}_2$ be cross-$t$-intersecting sub-families of $\mathcal{P}_{n,n}$. If $n$ is sufficiently large, then
\[|\mathcal{A}_1||\mathcal{A}_2| \leq ((n-t)!)^2,\]
and equality holds if $\mathcal{A}_1 = \mathcal{A}_2 = \left\{A
\in \mathcal{P}_{n,n} \colon \{(1,1), ..., (t,t)\} \subset A
\right\}$.
\end{theorem}
%
%\medskip \medskip
%\textbf{Acknowledgement.} The author is indebted to the anonymous
%referee for checking the paper carefully and providing remarks
%that led to an improvement in the presentation.

\end{document}